\documentclass{article}
\usepackage{tikz}
\usepackage{amssymb}
\usepackage{times}
\usepackage[T1]{fontenc}
\usepackage{pgf}
\usepackage{latexsym,slashbox}
\usepackage{color,xcolor}
\usepackage{graphicx}
\usepackage[latin1]{inputenc}
\usepackage{geometry} % Ò³ÃæÉèÖÃ
\usepackage{multirow}
\usepackage{enumerate}
\usepackage{amsmath}
\usepackage{manfnt}
\usepackage[english]{babel}
\usepackage{hyperref}
\usepackage{pifont}
\usepackage{amsmath,amssymb,amsfonts}
\newtheorem{prop}{Proposition}
\newtheorem{lem}{Lemma}
\newtheorem{thm}{Theorem}
\newtheorem{cor}{Corollary}

\newcommand{\F}{\mathbb{F}}
\newcommand{\C}{\mathcal C}

\begin{document}
\title{On error distance of received words with fixed degrees to Reed-Solomon code\footnote{A preliminary version can be seen in the first author's PhD thesis in 2008.}}
\author{Li Yujuan\\Science and Technology on Information\\ Assurance Laboratory\\Beijing,\   P.R.China \\liyj@amss.ac.cn\and Zhu Guizhen\\Data Communication Science and
\\Technology Research Institute\\Beijing,\   P. R.China\\zhugz08@gmail.com}
\date{}
\maketitle
\begin{abstract}
Under polynomial time reduction, the maximum likelihood decoding
of a linear code is equivalent to computing the error distance of a received word.
It is known that the decoding complexity of standard Reed-Solomon codes at certain
radius is at least as hard as the discrete logarithm problem over certain large
finite fields. This implies that computing the error distance is hard for standard
Reed-Solomon codes. Using some elegant algebraic constructions, we are able to determine the error distance of received words whose degree is $k+1$
to the Standard Reed-Solomon code or Primitive Reed-Solomon code exactly. Moreover, we can precisely determine the error distance of  received words of degree $k+2$ to the Standard Reed-Solomon codes. As a corollary, we can simply get the results of Zhang-Fu-Liao and Wu-Hong on the deep hole problem of Reed-Solomon codes.
%(add:over prime finite fields).
\end{abstract}
%\tableofcontents
\section{Introduction}
%There is always a possibility that a signal is corrupted when transferred over a long
%distance. Error-detecting and error-correcting codes alleviate the problem and make the
%modern communication possible. The Reed-Solomon codes are very popular in engineering
%(a reliable channel) due to their simplicity, burst error correction capabilities, and
%the powerful decoding algorithms within small error distance they admit.

Let~$\mathbb{F}_q$~be the finite field with~$q$~elements, where~$q$~is a prime power. For positive
integers~$k < n\leq q$, the generalized Reed-Solomon code, denoted by $\mathcal{C}$, can be thought
of as a map from~$\mathbb{F}^k_q \longrightarrow \mathbb{F}^n_q$, in which a message~$(a_0, a_1, \dots, a_{k-1})$~is mapped to a vector
~$( f (x_1), f (x_2), \dots , f (x_n))$, where~$f(x)=a_{k-1}x^{k-1} + a_{k-2}x^{k-2}+\dots+ a_0\in\mathbb{F}_q[x]$~and
$D=\{x_1, x_2,\dots , x_n\}\subseteq\mathbb{F}_q$~ is called the evaluation set. It is obvious that~$\mathcal{C}$~is a linear subspace
of~$\mathbb{F}^n_q$~with dimension~$k$.
 When the evaluation set is the whole field~$\mathbb{F}_q$, the resulting code is called the standard Reed-Solomon code, denoted by~$\mathcal{C}_q$. If the evaluation set is $\mathbb{F}_q^*$, the resulting code is called primitive Reed-Solomon code, denoted by $\mathcal{C}_q^*$.
 %In the literature, the standard
%Reed-Solomon code is often called the extended Reed-Solomon code. This name can
%be confused (to us) to the generalized Reed-Solomon code.

The Hamming distance between two codewords is the number of coordinates in
which they differ. The error distance of a received word~$u\in\mathbb{F}^n_q$~to the code~$\mathcal{C}$~is the
minimum Hamming distance of~$u$~to codewords, denoted by~$d(u, \mathcal{C})$. A Hamming ball of
radius $m$~is the set of vectors within Hamming distance $m$ to some vector in $\mathbb{F}^n_q$. The minimum
distance of a code is the smallest distance between any two distinct codewords,
and is a measure of how many errors the code can correct or detect. The covering radius
of a code is the maximum possible distance from any vector in $\mathbb{F}^n_q$
to the closest codeword.
A deep hole is a vector which achieves this maximum. The minimum distance
of generalized Reed-Solomon codes is $n-k+1$. The covering radius of generalized
Reed-Solomon codes is $n-k$. Therefore, all the deep holes of Reed-Solomon code are the vectors of error distance $n-k$.
\subsection{Related Work}
%The pursuit of efficient decoding algorithms for Reed-Solomon codes has yielded intriguing
%results. If the radius of a Hamming ball centered at a received word is less than
%half the minimum distance, there can be at most one codeword in the Hamming ball.
%Finding this unique codeword (if it exists) is called unambiguous decoding. It can be
%efficiently solved, see \cite{BW} for a simple algorithm. If the radius is somewhat larger, but
%less than $n-\sqrt{n(k-1)}$, the number of codewords is of polynomial size. In this case,
%the decoding problem can be efficiently solved by the Guruswami-Sudan list decoding
%algorithm \cite{GS}, which outputs all the codewords inside a Hamming ball. If the radius is
%stretched further, the number of codewords in a Hamming ball may be exponential.We
%then study the bounded distance decoding problem, which outputs just one codeword in
%a Hamming ball of certain radius. More importantly, we can remove the restriction on
%radius and investigate the maximum likelihood decoding problem, which is the problem
%of computing a closest codeword to a given vector in $\mathbb{F}^n_q$.

The complexity for decoding Reed-Solomon codes has also attracted attention recently.
Guruswami and Vardy \cite{GV} proved that the maximum likelihood decoding of
generalized Reed-Solomon codes is NP-hard. In fact, the weaker problem of deciding
deep holes for generalized Reed-Solomon codes is already co-NP-complete, see \cite{CM}.
In the much more interesting case of standard Reed-Solomon codes, it is unknown if
decoding remains NP-hard. This is still an open problem. Cheng and Wan \cite{CW} \cite{CW1} managed
to prove that the decoding problem of standard Reed-Solomon codes at certain
radius is at least as hard as the discrete logarithm problem over a large extension of a
finite field. This is the only complexity result that is known for decoding the standard
Reed-Solomon code.

Under polynomial time reduction, the maximum likelihood decoding
of a linear code is equivalent to computing the error distance of a received word.
Our aim of this paper is to study the problem of computing the error distance of received words of certain degrees to the Reed-Solomon code. We shall use algebraic methods. For this purpose, we
first define the notion of the degree of a received word.
For $u=(u_1,u_2,\dots,u_n)\in\mathbb{F}^n_q$, $D=\{x_1,\dots,x_n\}\subset \mathbb{F}_q$, let
\[ u(x) =\sum_{i=1}^nu_i\frac{\prod_{j\neq i}(x-x_j)}{\prod_{j\neq i}(x_i-x_j)}\in\mathbb{F}_q[x].\]

That is, $u(x)$ is the unique (Lagrange interpolation) polynomial of degree at most $n-1$
such that $u(x_i)=u_i$ for $1\leq i\leq n$. For $u\in\mathbb{F}^n_q$, we define $\deg(u) = \deg(u(x))$, called the
degree of $u$. It is clear that $d(u,\mathcal{C})=0$ iff $\deg(u)\leq k-1$. Without loss of generality, we can assume that
 $k\leq\deg(u)\leq n-1$ and $u(x)$ is monic. We have the following simple
bound.
\begin{lem}\label{lem1}
For $k\leq\deg(u)\leq n-1$, we have the inequality
\[n-\deg(u)\leq d(u, \C)\leq n-k.\]
\end{lem}

This result shows that if $\deg(u) = k$, then $d(u, \C) = n-k$ and thus $u$ is a deep hole.
As mentioned before, it is NP-hard to determine whether $d(u,\C) = n-k$ (the deep hole
problem) for generalized Reed-Solomon codes. Thus, one way to exploit this problem is restricting our attention to
the most natural and important case, namely the standard Reed-Solomon code $\C_q$. Even
in this restricted case, we cannot expect a complete solution to the problem of computing
the error distance, as it is at least as hard as the discrete logarithm in a large finite
field. However, we expect that a lot more can be said for standard Reed-Solomon codes.
For instance, Cheng and Murray \cite{CM} conjectured the following complete classification
of deep holes for standard Reed-Solomon codes.

\noindent
\textbf{Conjecture} (Cheng-Murray). All deep holes for standard Reed-Solomon codes are
those words satisfying $\deg(u) = k$. In other words, a received word u is a deep hole for
$\C_q$ iff $\deg(u) = k$.

The deep hole problem for generalized Reed-Solomon codes is NP-hard. In contrast,
the Cheng-Murray conjecture implies that the deep hole problem for the standard Reed-
Solomon code can be solved in polynomial time. A complete proof of this conjecture
(if correct) seems rather difficult at present. As a theoretical evidence, they proved that
their conjecture is true if $d := \deg(u)-k$ is small and $q$ is sufficiently large compared
to $d+k$. More precisely, they showed
\begin{prop}
 Let $u\in\F^q_q$ such that $1\leq d:= \deg(u)-k\leq q-1-k$. Assume that
$q\geq \max\{k^{7+\epsilon}, d^{\frac{13}{3}+\epsilon}\}$ for some constant $\epsilon> 0$. Then $d(u, \C_q)< q-k$, that is, $u$ is not a
deep hole.
\end{prop}
%The method of Cheng-Murray is to reduce the problem to the existence of a rational
%point on a hypersurface over $\F_q$. They showed that the resulting hypersurface is
%absolutely irreducible and then applied an explicit version of the Lang-Weil theorem.
However, they did not obtain the exact value of $d(u, \C_q)$, only the weaker inequality
$d(u, \C_q) < q-k$. Li and Wan \cite{LW} improved their results using Weil's character sum
estimate and the approach of Cheng-Wan \cite{CW} as follows.
\begin{prop}Let $u\in\F^q_q$ such that $1\leq d:= \deg(u)-k\leq q-1-k$. For some constant $\epsilon> 0$,
\begin{enumerate}[1)]
\item if
\[q\geq \max\{(k+1)^2, d^{2+\epsilon}\},\text{ and }k>\left(\frac{2}{\epsilon}+1\right)d+\frac{8}{\epsilon}+2,\]
then $u$ is not a deep hole.
\item if
\[q\geq \max\{(k+1)^2, (d-1)^{2+\epsilon}\},\text{ and }k>\left(\frac{4}{\epsilon}+1\right)d+\frac{4}{\epsilon}+2,\]
then $d(u,\C_q)=q-(k+d).$
\end{enumerate}
\end{prop}

Note that the last part of the proposition determines the exact error distance $d(u, \C_q)$
under a suitable hypothesis. Using a similar character sum approach, Qunying Liao
\cite{LIAO} unified the above two results of Li-Wan and proved the following extension.
\begin{prop}
Let $r\geq 1$ be an integer. For any received word $u\in\F^q_q$,$r\leq d:= \deg(u)-k\leq q-1-k$.
If
\[q\geq \max\left\{2{k+r\choose 2}+d, d^{2+\epsilon}\right\},\text{ and }k>\left(\frac{2}{\epsilon}+1\right)d+\frac{4+2r}{\epsilon}+2,\]
for some constant $\epsilon > 0$, then $d(u, \C_q)\leq q-k-r$.
\end{prop}

Antonio Cafure etc. \cite{ANT} uses a much more sophisticated algebraic geometry
approach and obtains a slightly improvement of one of the Li-Wan results.
\begin{prop}
Let $u\in\F^q_q$ such that $1\leq d:= \deg(u)-k\leq q-1-k$. Assume that
\[q\geq \max\{(k+1)^2, 14d^{2+\epsilon}\},\text{ and, }k>\left(\frac{2}{\epsilon}+1\right)d,\]
for some constant $\epsilon> 0$, then $u$ is not a deep hole.
\end{prop}

Again, this result gives only the inequality $d(u, \C_q)<q-k$, not the exact value of the
error distance $d(u, \C_q)$.
As for the error distance, Zhu-Wan \cite{ZW2} prove the following result.
 %Using the Weil bound and a new sieve for distinct coordinates
%counting, they are able to compute the error distance for a large class of
%received words.
\begin{prop}
 Let $r\geq 1$ be an integer and $u\in\F^q_q$, $r\leq d:= \deg(u)-k\leq q-1-k$. There are positive constants $c_1$ and $c_2$ such that if
 \[d<c_1 q^{1/2}, \left(\frac{d+r}{2}+1\right)\log_2q<k<c_2q,\]
 then $d(u, \C_q)\leq q-k-r$.
\end{prop}
So far, no one has proved or defied Cheng-Murray's conjecture on standard Reed-Solomon code.  In a recent paper by Cheng-Li-Zhuang\cite{CLZ}, they classify deep holes
completely for generalized Reed-Solomon codes $\mathcal{C}$, where $q$ is
prime, and $D>k\geq\frac{q-1}{2}$. Moreover, they prove that
\begin{prop}
Cheng-Murray's conjecture is true for $p>2, k+1\leq p $ or $3\leq q-p+1\leq k+1\leq q-2$.
\end{prop}

Another way to research the maximum likely decoding problem or deep hole problem is studying the error distance of received words of certain degree to Reed-Solomon code. Wu-Hong\cite{WH} show that some received words of degree $q-2$ are deep holes of Primitive Reed-Solomon code.
\begin{prop}
Let $\C^*_q$ be a Primitive Reed-Solomon code, $q\geq 4$, and $2\leq k\leq q-2$, then all the received words that can be presented by polynomial $u(x)=ax^{q-2}+v(x)$ are deep holes, where $\deg(v)\leq k-1$, $a\neq 0$.
\end{prop}
Zhang-Fu-Liao\cite{ZFL} extend the result above to any evaluation set $D\neq \F_q$, and derive the following conclusion.
\begin{prop}
Let $\C$ be a Generalized Reed-Solomon code, $D\neq \F_q$, then for $a\neq 0$, $b\notin D$, all the received words that can be presented by polynomial $u(x)=a(x-b)^{q-2}+v(x)$ are deep holes, where $\deg(v)\leq k-1$.
\end{prop}
Meanwhile, they find another kind of deep holes for certain dimension and finite fields.
\begin{prop}
Let $q$ be a power of 2, and $q\geq4$, $\C$ be a Generalized Reed-Solomon code, evaluation set $D=\F_q^*$ or $D=\F_q^*\setminus \{1\}$, $k=q-4$. If $a\neq 0$, then  all the received words that can be presented by polynomial $u(x)=ax^{q-3}+v(x)$ are deep holes, where $\deg(v)\leq k-1$.
\end{prop}
Finally, they prove that if $q>5$ with odd characteristic $p$, and $2\leq k\leq q-3$, all the received words represented by the following polynomials are not deep holes of Primitive Reed-Solomon code.
\[u(x)=ax^{k+2}+bx^{k+1}+cx^k+v(x),\]
where $a\in\F_q^*$, $b,c\in\F_q$, and $\deg(v)\leq k-1$.

Again, this result gives only the inequality $d(u, \C_q)<q-k$, not the exact value of the
error distance $d(u, \C_q)$. And they only discuss the case that characteristic $p\neq 2$. For $p=2$, the problem may be more complicated, which can be seen in our analysis in this paper.
%In fact,  for $p=2$, we find a new deep hole of degree $k+2$. That is for $q=8$, $k=1$ or $k=2$, the received word of polynomial $x^{k+2}+v(x)$ with $\deg(v(x))\leq k-1$ is a deep hole of Standard Reed-Solomon code. This is a new kind of deep hole and Zhang-Fu-Liao miss it in \cite{}. We will give it in 1.2.

\subsection{Our results}

 In this paper, we focus on computing the error distance of received words of fixed degrees to Reed-Solomon codes. The main results consist of two parts. Firstly, we exploit the error distance of received words of degree $k+1$ to Standard Reed-Solomon code and Generalized Reed-Solomon code; Secondly, we compute the error distance of received words of degree $k+2$ to Standard Reed-Solomon code not only for $p\neq 2$ but also for the case $p=2$. As a corollary, we can rather easily get Wu-Hong and Zhang-Fu-Liao's results on deep hole.
\begin{thm}\label{thm1}
\begin{enumerate}[(i)]
\item Let $\C_q$ be a Standard Reed-Solomon code, and $u\in\F_q^n$ represented by polynomial $u(x)=x^{k+1}-bx^k+v(x)$, $\deg(v)\leq k-1$, then  $d(u,\C_q)=q-k$ if one of the following holds
    \begin{enumerate}
    \item $b=0, p=2, k=1$,
    \item  $b=0, p=2, k=q-3$,
    \end{enumerate}
    otherwise, $d(u,\C_q)=q-k-1$.
  % \begin{align*}
%   d(u,\C_q)=\left\{\begin{array}{ll}
%    q-k, & b=0, p=2, k=1 \text{ or } q-3\\
%    q-k-1, & otherwise.
%    \end{array}\right.
%   \end{align*}
\item If $D=\F_q^*$  and $q>5$, then $ d(u,\C)=q-k-1$ if one of the following holds
\begin{enumerate}
\item $b=0, p=2, k=1 $,
\item $b=0, p=2, k=q-4$,
\item $b=0, k=q-3$.
\end{enumerate}
otherwise, $ d(u,\C)=q-k-2$.
 %\begin{align*}
%   d(u,\C)=\left\{\begin{array}{ll}
%    q-k-1, & b=0, p=2, k=1 \text{ or } q-4 \text{ or } b=0, k=q-3\\
%    q-k-2, & otherwise.
%    \end{array}\right.
%   \end{align*}
\end{enumerate}
\end{thm}

\begin{thm}
Let  $\C_q$ be a Standard Reed-Solomon code, $k\geq 1, k+2\leq q-1$, and $u\in\F_q^n$ represented by polynomial $u(x)=x^{k+2}-bx^{k+1}+cx^k+v(x)$, $\deg(v)\leq k-1$, then
\begin{enumerate}[(i)]
\item If $k+2=q-1$, then
\begin{align*}
d(u,\C_q)=\left\{
\begin{array}{ll}
q-k-2 &\text{if  } b^{2}=c,\\
q-k-1 &\text{if  } b^{2}\neq c.
\end{array}
\right.
\end{align*}

\item  If $p=2$ and $ k+2\leq q-2$, then we can get the following results.
\begin{enumerate}
\item $d(u,\C)=q-k-2$ if $(k,b,c)$ satisfies one of the following conditions.
  \begin{itemize}
    \item $2\mid k+1, 4\nmid k+1$, and $b^2\neq c $.
    \item $2\mid k+1, 4\nmid k+1$, $b^2= c$ and $k+2>q/2$.
    \item $4\mid k+1$ and $c\neq 0$.
   \item $4\mid k+1$,$c=0$ and $k+2<q/2$.
    \end{itemize}
\item $d(u,\C)\leq q-k-1$ if $(k,b,c)$ satisfies one of the following conditions.
\begin{itemize}
\item $4\mid k+1$, $c=0$ and $k+2\geq q/2$.
\item $4\mid k$.
\item $2\mid k$,$4\nmid k$, and $b\neq 0$.
\item $2\mid k$,$4\nmid k$, and  $c\neq 0$.
\item $2\mid k$,$4\nmid k$, $b=c=0$ and $k+1>q/2$.
\end{itemize}
\end{enumerate}

\item  If $p\neq 2$ and $ k+2\leq q-2$, then if $p\nmid k+2$, we have $d(u, \mathcal{C})\leq q-k-1.$ In the case that $p\mid k+2$ we can conclude the following results.
    \begin{enumerate}
     \item If $b=c=0$ and $k+2>q/2+1$, then $ \ d(u, \mathcal{C})\leq q-k-1$;
     \item If $b\neq 0$, then $ \ d(u, \mathcal{C})=q-k-2$;
     \item If $c\neq 0$, then
     \begin{align*}
\ d(u, \mathcal{C})=\left\{\begin{array}{lll}q-k &\text{if } p=3, k+2=3 \text{ and }-c \text{ is not a nonzero square},\\
q-k-1 & \text{if } p=3, k+2=q-3 \text{ and }-c  \text{ is not a nonzero square},\\
q-k-2  &otherwise.
\end{array}
\right.
\end{align*}
    \end{enumerate}

\end{enumerate}
\end{thm}

In particular, for the cases which do not satisfy the conditions we discuss in our theorem, we find some new deep holes.
\begin{itemize}
\item $q=8$, $k=1$, $b^2=c\in\mathbb{F}_8$. Received word with polynomial $u=x^3+bx^2+cx+d$ is a deep hole of $\mathcal{C}_8=\{(x_i,x_i,\dots,x_i)\in\mathbb{F}_8^8\mid x_i\in\mathbb{F}_8, 1\leq i\leq8\}$, where $d\in\mathbb{F}_8$.
\item $q=8$, $k=2$, $b^2=c=0$.Received word with polynomial $u=x^4+dx+e$ is a deep hole of $\mathcal{C}_8=\{(mx_1+t,mx_2+t,\dots,mx_8+t)\in\mathbb{F}_8^8\mid x_i\in\mathbb{F}_8, 1\leq i\leq8, m, t\in\mathbb{F}_8\}$, where $d,e\in\mathbb{F}_8$.
\end{itemize}
% By theorem 2, we can have the results of deep holes.
%\begin{cor}
%If q is a prime not equal to 2 and $k+1\leq q-1$, then $d(u,\C_q)\leq q-k-1$. In other words, all the received words with degree $k+2$ are not deep holes over odd prime finite fields.
%\end{cor}
In our proof, we convert the problem of deciding the error distance of a received word
to solving a polynomial equation. Compared with  approach in [2][3],
our method is much simpler and using some algebraic constructions and character sum estimate, we not only get the deep hole results, but also can determine the error distance explicitly.

\emph{Organization.} In Section 2, we provide a brief introduction to finite fields and state some fundamental definitions and lemmas. In Section 3, we discuss the error distance of received words of degree $k+1$ to Standard Reed-Solomon code and Primitive Reed-Solomon code respectively. The case that computing error distance of received words of degree $k+2$ to Standard Reed-Solomon code is studied in Section 4. %Finally, Section 5 comes into a conclusion.
\section{Preliminaries}
We first review the theory of finite field and character sums in the form we need. Let $\F_q$ be the finite field with character $p$, where $q$ is a $p$ power.
For a element $a\in\F_q^*$, the order of $a$ is defined by the smallest number $d$ such that $a^d=1$.
Let $\chi: \mathbb{F}_q^*\longrightarrow\mathbb{C}^*$ be a multiplicative character from the invertible elements
of $\F_q$ to the non-zero complex numbers and satisfies that for $a,b\in \mathbb{F}_q^*$, $\chi(ab)=\chi(a)\chi(b)$. If for all $a\in\mathbb{F}_q^*$, $\chi(a)=1$, then call $\chi$ trivial character, denoted by 1. The smallest $d$ such that $\chi^d=1$ is called the degree of $\chi$. Extend the definition to $\mathbb{F}_q$ by
\[\chi(0)=\left\{\begin{array}{ll}1,&\mbox{$\chi=1$};\\
0,&\mbox{$\chi\neq 1$}\end{array}\right. \]
\begin{lem} \label{2.1}\cite{LH} Let $\mathbb{F}_{q}$ be a finite field, $p\neq
2$. If $n$ is odd and $a_i\neq
0, 1\leq i\leq n$, then the number of solutions of equation $a_1x_1^2+\cdots+a_nx_n^2=b$ over $\mathbb{F}_{q}$ is
$$q^{n-1}+q^{(n-1)/2}\eta((-1)^{(n-1)/2}ba_1\cdots a_n),$$
where $\eta$ is a character of degree 2 over $\mathbb{F}_{q}$.
\end{lem}
\begin{lem}\label{2.2} \cite{LH}Let $\mathbb{F}_{q}$ be a finite field, $p\neq
2$. If $n$ is even and $a_i\neq
0, 1\leq i\leq n$, then the number of solutions of equation $a_1x_1^2+\cdots+a_nx_n^2=b$ over $\mathbb{F}_{q}$ is
$$q^{n-1}+v(b)q^{(n-2)/2}\eta((-1)^{n/2}a_1\cdots a_n),$$
where $\eta$ is a character of degree 2 over $\mathbb{F}_{q}$.
\end{lem}
%ÏÂÃæµÄÒýÀíÊÇÎÒÃÇÔÚÖ¤Ã÷ÖÐ¾­³£ÒªÓÃµ½µÄ¡£¼òÊöÈçÏÂ£º
\begin{lem} \label{2.3}\cite{LW}
Let
$u\in \mathbb{F}_{q}^n$ be a received word with degree $k+r$, where $k+1\leq k+r\leq n-1 $.
Then
\begin{enumerate}[(i)]
\item $d(u, \mathcal{C})=n-k-r
$ if and only if there exists a subset $E=\{x_{1},\dots,
x_{k+r}\}$ of $\mathcal{D}$ such that
$$u(x)-v(x)=(x-x_{1})\cdots(x-x_{k+r}),$$
for some $v(x)\in\mathbb{F}_{q}[x], \ \deg v(x)\leq k-1$.
\item $d(u, \mathcal{C})\leq n-k-i
,\ (1\leq i\leq r)$ if and only if there exists a subset $E=\{x_{1},\dots,
x_{k+i}\}$ of $\mathcal{D}$ and a monic polynomial$g(x)$ of degree $r-i$ such that $$u(x)-v(x)=(x-x_{1})\cdots(x-x_{k+i})g(x),$$
for some $v(x)\in\mathbb{F}_{q}[x], \ \deg v(x)\leq k-1$.
\end{enumerate}

\end{lem}
\section{The case for received words of degree $k+1$}
In this section, we give the proof of Theorem \ref{thm1}. Let $u\in\F_q^n$  be a received word represented by polynomial $u(x)=x^{k+1}-bx^k+v(x)$ with $\deg(v)\leq k-1$.
\subsection{Computing $d(u,\C_q)$}
From Lemma \ref{lem1} and  Lemma 4, we know that $q-k-1\leq d(u,\C_q)\leq q-k$, and $d(u,\C_q)=q-k-1$ if and only if
there exists a subset $\{x_1, x_2,\dots, x_{k+1}\}\subset\F_q$ of size $k+1$ such
that $b=x_1+\dots+x_{k+1}$.
\begin{itemize}
\item $b\neq 0$.

Let $g$ be a primitive element in $\F_q$, and
\[u=1+g+g^2+\dots+g^{k-1}+0.\]
as the order of $g$ is  $q-1$ and $k+1\leq q-1$, then $u\neq 0$ and
\[1=0+u^{-1}+u^{-1}g+\dots+u^{-1}g^{k-1},\]
thus,
\[b=0+bu^{-1}+bu^{-1}g+\dots+bu^{-1}g^{k-1},\]
obviously, the $k+1$ items above are distinct. To be concluded, if $b\neq 0$, $d(u,\C_q)=q-k-1$.
\item $b=0$.
\begin{itemize}
\item $b=0, p\neq 2$.%(??If$p\neq 2$)

In this case, for any $x\in\F_q^*$, $x\neq -x$, therefore
\[\F_q=\{0,x_1, -x_1, \dots, x_{\frac{q-1}{2}}, -x_{\frac{q-1}{2}}\}.\]
If $k$ is odd, then
\[0=x_1+(-x_1)+\dots+x_{\frac{k+1}{2}}+(-x_{\frac{k+1}{2}}).\]
If $k$ is even, we only need to add 0 to the right side of the equation above. Therefore, if $b=0, p\neq 2$, $d(u,\C_q)=q-k-1$.

\item $b=0, p=2$.%(??If$p=2$)

Without loss of generality, we can assume $q>2$. As the sum of all the elements in $\F_q$ is 0,  the conclusion holds for $k+1$ iff  it holds for $q-k-1$. For $k+1=2$, $d(u,\C_q)=q-k-1$ is equivalent to the fact that there exists $x_1, x_2\in\F_q$, $x_1\neq x_2$ and $x_1+x_2=0$. But when $p=2$, if $x_1+x_2=0$, then $x_1=x_2$, contradiction. Thus, if $p=2, b=0$, and $k+1=2$, $d(u,\C_q)=q-k$. Likewise, if $p=2, b=0$, and $q-k-1=2$, $d(u,\C_q)=q-k$. Without loss of generality, we can assume $2<k+1\leq q/2$. Set
\[S=\F_q^*\setminus \{g, g^2, \dots, g^{k-1}\}.\]
It is easy to see that the number of elements in $S$ is $q-1-(k-1)=q-k>q/2$. Set
\[T=\{g+g^2+\dots+g^{k-1}+g^i\mid g^i\in S\},\]
then the number of elements in $T$ is also $q-k$, thus, $|S|+|T|>q$, and $S\cup T\subseteq\mathbb{F}_q$, which means that there exist two elements $g^i$ and $g^j$ in $S$ such that
\[g+g^2+\dots+g^{k-1}+g^i=g^j.\]
As $p=2$, then we have
\[g+g^2+\dots+g^{k-1}+g^i+g^j=0.\]
Obviously, these $k+1$ elements are distinct, so far, we can conclude that $d(u, \C_q)=q-k-1$. Likewise, the same conclusion holds for $2<q-k-1\leq q/2$. Overall, if $p=2, b=0$ and $1<k<q-3$, then $d(u, \C_q)=q-k-1$. As for the case $k=q-2$, for any $b$, the sum of  $q-1$ elements in $\F_q\setminus\{-b\}$ is $b$. Thus, if $k=q-2$, $d(u,\C_q)=q-k-1$. The proof of the first part of Theorem \ref{thm1} is complete.
\end{itemize}
\end{itemize}
\subsection{Computing $d(u,\C_q^*)$ }
The proof of Theorem \ref{thm1}(ii) is similar to the proof of Theorem \ref{thm1}(i). %For simpleness(??So), we only give the proof for $D=\F_q^*$.
\begin{itemize}
\item $b\neq0$.
For this case, the proof is the same as the proof for $b\neq0$ in section 3.1. We omit it and conclude that if $b\neq 0$, $d(u,\C_q^*)=q-k-2$.

%Let $g$ be a primitive element in $\F_q$, and
%\[u=1+g+g^2+\dots+g^{k-1}+g^k.\]
%As the order of $g$ is  $q-1$ and $k+1\leq q-2$, then $u\neq 0$ and
%\[1=u^{-1}+u^{-1}g+\dots+u^{-1}g^{k-1}+u^{-1}g^k,\]
%thus,
%$$b=bu^{-1}+bu^{-1}g+\cdots+bu^{-1}g^{k},$$
%obviously, these $k+1$ items above are distinct. To be concluded, if $b\neq 0$, $d(u,\C_q)=q-k-2$.

\item $b=0$.
\begin{itemize}
\item $b=0$, $k+1=q-2$.

If there exist $k+1$ distinct elements $x_1,x_2,\dots,x_{k+1}$ in $\F_q^*$ satisfying $\sum x_i=0$, there is only one nonzero element in $\F_q^*\setminus\{x_1,x_2,\dots,x_{k+1}\}$, which contradicts the fact that the sum of all elements in $\F_q^*$ is 0. Thus, if $b=0$ and $k+1=q-2$, $d(u,\mathcal{C})=q-k-1$.
\item $b=0$, $k+1\neq q-2$.

If $p\neq2$, and $k+1$ is even, the proof is same as the proof for the case $b=0, p\neq2$, and $k+1$ is even in Section 3.1. Now we discuss the case that $k+1$ is odd. As $q>5$, we can find  $z_{1},z_{2}\in\mathbb{F}_{q}^{\ast}$ satisfying $z_{1}\neq z_{2},
-z_{2},-2z_{2},-\frac{1}{2}z_{2}$. Thus,
$$\mathbb{F}_{q}=\{z_{1}+z_{2}+z|z\in\mathbb{F}_{q}\}.$$
So there exists $z_{3}\in\mathbb{F}_{q}^{\ast}$ such that $z_{1},z_{2},z_{3}$ are distinct and
$$z_{1}+z_{2}+z_{3}=0.$$
As the sum of all the elements in $\F_q^*$ is 0,  the conclusion holds for $k+1$ iff  it holds for $q-k-2$.
As $k+ 1$ is odd and $k+1\neq q-2$, we can assume $k+1\leq q-4$, say $k-2\leq q-7$. Set
$$M=\mathbb{F}_{q}^{\ast}\setminus\{\pm z_{1},\pm z_{2},\pm z_{3}\}$$
then the number of elements in $M$ is also $q-7$, together with the fact that $k-2$ is even, we can get $k-2$ elements in $M$ summing to 0 similarly to what we proved in last subsection. Adding $z_{1},z_{2},z_{3}$ into  these $k-2$ elements, then we get $k+1$ distinct elements in $\F_q^*$ whose sum is 0. If $p=2$, the proof is same as the proof for the case $p=2,b=0$ in section 3.1. %Using the similar proof as Section 3.1, we can conclude if $k+1=2$ or $q-1-(k+1)=2$(i.e. $k=q-4$), $d(u,\mathcal{C})=q-k-1$.
\end{itemize}
\end{itemize}
This completes the proof of Theorem \ref{thm1}.

\section{The case for received words of degree $k+2$}

\begin{lem}\label{lem2}
\begin{enumerate}[(i)]
\item Suppose $p=2$, $2\leq t\leq q-2$, and $c\in\mathbb{F}_{q}^{\ast}$. Then
there exist $t$ distinct elements $\gamma_{1},\gamma_{2},\ldots,\gamma_{t}$ in $\mathbb{F}_{q}^{\ast}$ such that
$$c=\sum_{1\leq
i<j\leq t}\gamma_{i}\gamma_{j}.$$
\item Suppose  $p=2$, $2\leq t\leq q-3$, and $c\in\mathbb{F}_{q}^{\ast}$,
then there exist $t$ distinct elements $\gamma_{1},\gamma_{2},\ldots,\gamma_{t}$ in $\mathbb{F}_{q}^{\ast}$ such that $$c=\sum_{1\leq
i\leq j\leq t}\gamma_{i}\gamma_{j}.$$
\item If ~$t=q-1$, then $\sum\limits_{1\leq i<j\leq q-1\atop
\gamma_{i},\gamma_{j}\in\mathbb{F}_{q}^{\ast}}\gamma_{i}\gamma_{j}=0$.
\end{enumerate}
\end{lem}
\emph{Proof.} Let $g$ be a primitive element in $\mathbb{F}_{q}$.
\begin{enumerate}[(i)]
\item  As $2\leq t\leq q-2$,
$(1-g^{t-1})(1-g^{t})\neq 0$. Therefore, for $p=2$ and any $c\in\mathbb{F}_{q}^{\ast}$,
the following equation with variable $y$ always has solutions.
\begin{equation*}
y^{2}\frac{1}{1-g}\frac{g(1-g^{t-1})(1-g^{t})}{1-g^{2}}=c.
\end{equation*}
Suppose $\sigma$ is a root of the equation above
and set $\gamma_{i}=\sigma g^{i-1},1\leq i\leq
t$, then
\begin{eqnarray*}
\sum_{1\leq i<j\leq t}\gamma_{i}\gamma_{j}
&=&\sigma^{2}\sum_{0\leq i<j\leq
t-1}g^{i+j}\\
&=&\sigma^{2}[(g+g^{2}+\cdots+g^{t-1})+g(g^{2}+\cdots+g^{t-1})+\cdots\\&+&g^{t-4}(g^{t-3}+g^{t-2}+g^{t-1})+g^{t-3}(g^{t-2}+g^{t-1})+g^{t-2}g^{t-1}]\\
%&=&\sigma^{2}\frac{1}{1-g}\{g-g^{t}+g^{3}-g^{t+1}+\cdots\\&+&g^{2t-7}-g^{2t-4}+g^{2t-5}-g^{2t-3}+g^{2t-3}-g^{2t-2}\}\\
&=&\sigma^{2}\frac{1}{1-g}\{(g+g^{3}+\cdots+g^{2t-7}+g^{2t-5}+g^{2t-3})\\&-&(g^{t}+g^{t+1}+\cdots+g^{2t-4}+g^{2t-3}+g^{2t-2})\}\\
&=&\sigma^{2}\frac{1}{1-g}\frac{g(1-g^{t-1})(1-g^{t})}{1-g^{2}} \\
&=&c.
\end{eqnarray*}
\item According to $2\leq t\leq q-3$, we can conclude that $(1-g^{t+1})(1-g^{t})\neq
0$. Then for any $c\in\mathbb{F}_{q}^{\ast}$,
let $\omega$ be a solution of the following equation.
\begin{equation*}
y^{2}\frac{1}{1-g}\frac{(1-g^{t+1})(1-g^{t})}{1-g^{2}}=c.
\end{equation*}
Set $\gamma_{i}=\omega g^{i-1},1\leq i\leq t$. Similarly to (i), we can deduce that
\begin{eqnarray*}
\sum_{1\leq i<j\leq t}\gamma_{i}\gamma_{j}
&=&\omega^{2}\sum_{0\leq i\leq j\leq
t-1}g^{i}g^{j}\\
&=&\omega^{2}\{1+g^{2}+\cdots+g^{2(t-1)}+\sum_{0\leq i<j\leq
t-1}g^{i+j}\}\\
&=&\omega^{2}\frac{1}{1-g}\frac{(1-g^{t+1})(1-g^{t})}{1-g^{2}}\\
&=&c.
\end{eqnarray*}
\item For the case that $t=q-1$, all the elements in $\F_q^*$ have the form $g^i, \ i=1,2,\dots, q-2$, so

\begin{eqnarray*}
\sum_{0\leq i<j\leq
q-2}g^{i}g^{j}&=&\sum_{0\leq i<j\leq
q-2}g^{i+j}\\
&=&\frac{1}{1-g}\frac{g(1-g^{q-2})(1-g^{q-1})}{1-g^{2}}\\
&=&0
\end{eqnarray*}
%\qed
%
\end{enumerate}
\begin{lem}\label{lem3}
Assume that $1\leq t<\frac{q}{2}-1$. If $p=2$, suppose $4\mid t$ and if $p\neq
2$, suppose $p\mid t$, then there exist $t$ distinct elements $\xi_{1},\dots,
\xi_{t}$ in $\mathbb{F}_{q}^{\ast}$ such that $$\sum_{1\leq i<j\leq t}\xi_{i}\xi_{j}=0.$$
\end{lem}
\emph{Proof.} Let $g$ be a primitive element in $\mathbb{F}_{q}$, set
$$M_{1}=\frac{1-g^{t-1}}{1-g},
M_{2}=\frac{1}{1-g}\frac{g(1-g^{t-2})(1-g^{t-1})}{1-g^{2}}.$$
 Let $\phi\in\mathbb{F}_{q}^{\ast}$ satisfy $\phi\neq
-M_{1},\phi^{2}+2\phi M_{1}+M_{2}\neq 0$,  for $2\leq i\leq t$, $\phi\neq g^{i-2}$ and
if $g^{i-2}+M_{1}\neq 0,$ $$\phi\neq
-\frac{M_{2}+g^{i-2}M_{1}}{g^{i-2}+M_{1}}.$$
As $1\leq t< \frac{q}{2}-1$,
$$(q-1)-1-2-2(t-1)=q-2-2t>0.$$
So $\phi$ does exist. Set
\begin{equation*}
z=\frac{\phi M_{1}+M_{2}}{\phi+M_{1}}
\end{equation*}
and $\xi_{1}=z+\phi,\ \xi_{i}=z+g^{i-2}, 2\leq i\leq t.$
Because of the choice of $\phi$, we can easily varify that $\xi_{1} \neq 0$, and $\xi_i\neq\xi_j$, $i\neq j$.
For $2\leq i\leq t$, if $g^{i-2}+M_{1}\neq 0$, obviously, $z+g^{i-2}\neq 0$,
i.e. $\xi_{i}\neq 0$. If $g^{i-2}+M_{1}=0$,
\begin{eqnarray*}
\xi_{i}&=&z+g^{i-2}\\
&=&\frac{\phi M_{1}+M_{2}}{\phi+M_{1}}+g^{i-2}\\
&=&\frac{\phi M_{1}+M_{2}}{\phi+M_{1}}-M_{1}\\
&=&\frac{M_2-M_1^2}{\phi+M_{1}}.
\end{eqnarray*}
We have
\begin{eqnarray*}
M_{2}-M_{1}^{2}
&=&\frac{1}{1-g}\frac{g(1-g^{t-2})(1-g^{t-1})}{1-g^{2}}-\frac{(1-g^{t-1})^{2}}{(1-g)^{2}}\\
&=&\frac{(1-g^{t-1})(g^{t}-1)}{(1-g)^{2}(1+g)}\\
&\neq&0.
\end{eqnarray*}
Then $\xi_i\neq 0$.
Under the condition that $p=2$, $4\mid t$, or $p\neq 2$, $p\mid t$, we can get
\begin{eqnarray*}
\sum_{1\leq i<j\leq t}\xi_{i}\xi_{j}&=&(z+\phi)\sum_{2\leq i\leq
t}(z+g^{i-2})+\sum_{2\leq i<j\leq
t}(z+g^{i-2})(z+g^{j-2})\\
&=&\frac{t(t-1)}{2}z^{2}+(t-1)(\phi+1+g+\cdots+g^{t-2})z+\phi
M_{1}+\sum_{2\leq
i<j\leq t}g^{i+j-4}\\
&=&-(\phi+M_{1})z+\phi M_{1}+
M_{2}\\
&=&0.
\end{eqnarray*}
%\begin{cor}
%Assume that $1\leq t<\frac{q}{2}$. If $p=2$, suppose $4\mid t-1$ and if $p\neq
%2$, suppose $p\mid t-1$, then there exist $t$ distinct elements $\xi_{1},\dots,
%\xi_{t}$ in $\mathbb{F}_{q}$ such that $$\sum_{1\leq i<j\leq t}\xi_{i}\xi_{j}=0.$$
%\end{cor}
By lemma 5 (iii) and lemma 6, we can easily get the following corollary.
\begin{cor}\label{cor1}
Assume that $ t>\frac{q}{2}$. If $p=2$, suppose $4\mid q-1-t$ and if $p\neq
2$, suppose $p\mid q-1-t$, then there exist $t$ distinct elements $\xi_{1},\dots,
\xi_{t}$ in $\mathbb{F}_{q}^*$ such that $$\sum_{1\leq i\leq j\leq t}\xi_{i}\xi_{j}=0.$$
\end{cor}
%\qed
%
\begin{lem}\label{4.3}
Assume that $p\neq 2$, $r,r_{1},\mu \neq 0$, and $b,c\in \mathbb{F}_{q}$, denote
$A=\big\{\frac{\alpha^{2}}{\mu}-\frac{b^{2}r^{2}}{\mu}+2cr_{1}~|~
\alpha \in \mathbb{F}_{q}\big\}$.
If $2<t<\frac{q+1}{2}$ and $t$ is even, then there exist $t$ distinct elements $y_{1}, \cdots,
y_{t}$ in $\mathbb{F}_{q}^*$ such that
\begin{equation*}(y_{1}+\cdots+
y_{t})^{2}-r(y_1^2+\cdots+y_{t}^2)\in A.
\end{equation*}
\end{lem}
\emph{Proof.} Let $g$ be a primitive element of $\mathbb{F}_{q}$. From Lemma \ref{2.1}, the following equation with variables
$\alpha,y,z$ has at most $(q^{2}-1)$  nonzero solutions in $\mathbb{F}_{q}^3$.
\begin{equation}\label{lem6}
\alpha^{2}+g^{2}z^{2}\frac{1-g^{t-2}}{1-g^{2}}-
g^{2}y^{2}\frac{1-g^{t}}{1-g^{2}}=0
\end{equation}
Denote
$$T=\{(\alpha,z)\in\mathbb{F}_{q}^2\mid \alpha,z\neq 0,
\text{ and there exists }y\in\mathbb{F}_{q}^{\ast}\text{ such that }(\alpha,z,y)\text{is a nonzero solution of Equation (\ref{lem6})}\}.$$
Thus, there are at most $\frac{q^2-1}{2}$ elements in $T$. The number of pairs $(\alpha,
z)\in\mathbb{F}_{q}^2$ such that $\alpha\neq 0$ and $\alpha\neq
\pm zg^i, 1\leq i\leq (t-2)/2, z\neq 0$ is $(q-1)(q-t+1)$. %\textcolor{red}{(ÎÒ¾õµÃÊÇ$(q-1)(q-t+1)$????)}
For $\frac{q+1}{2}>t>2$,  $$(q-1)(q-t+1)-\frac{q^{2}-1}{2}> 0,$$
then there exist $ \alpha_{1}, z_{1}\in\mathbb{F}_{q}^\ast, \alpha_{1}\neq
\pm z_{1}g^i, 1\leq i\leq (t-2)/2$, and for all $y\in
\mathbb{F}_{q}^*$
\begin{equation}{\alpha_{1}}^{2}(-2r)+(-2r)g^{2}z^{2}\frac{1-g^{t-2}}{1-g^{2}}\neq
(-2r)g^{2}y^{2}\frac{1-g^{t}}{1-g^{2}}.
\end{equation}
Set $y_{1}=yg, y_{2}=-yg, \dots, y_{t-1}=yg^{\frac{t}{2}},
y_{t}=-yg^{\frac{t}{2}}, y\in \mathbb{F}_{q}^*$, $m=y_1+y_2+\dots+y_t$, $n=y_1^2+y_2^2+\dots+y_t^2$, we have
$$m^{2}-rn=(-2r)g^{2}y^{2}\frac{1-g^{t}}{1-g^{2}}.$$
Set $y_{1}=\alpha_{1}, y_{2}=-\alpha_{1}, y_{3}=z_{1}g,
y_{4}=-z_{1}g,\dots, y_{t-1}=z_{1}g^{\frac{t-2}{2}},
y_{t}=-z_{1}g^{\frac{t-2}{2}}$,  $m=y_1+y_2+\dots+y_t$, $n=y_1^2+y_2^2+\dots+y_t^2$, then
$$m^{2}-rn=(-2r){\alpha_{1}}^{2}+(-2r)g^{2}{z_{1}}^{2}\frac{1-g^{t-2}}{1-g^{2}}.$$
According to (2), denote
$$B=\left\{(-2r)g^{2}y^{2}\frac{1-g^{t}}{1-g^{2}}\mid y\in
\mathbb{F}_{q}^*\right\}\cup
\left\{(-2r)({\alpha_{1}}^{2}+g^{2}{z_{1}}^{2}\frac{1-g^{t-2}}{1-g^{2}})\right\}.$$
Therefore, $|A|+|B|=q+1,$ which implies that $A\cap B\neq \emptyset$. In other words,
there exist $t$ distinct elements $y_{1}, \cdots,
y_{t}$ in $\mathbb{F}_{q}^*$ such that
$m^2-rn\in A.$

\subsection{Proof of theorem 2 (i).}
\begin{enumerate}[i)]
\item $b^2=c$.

From Lemma  \ref{2.3}, we know that $d(u, \mathcal{C})=q-k-2$ iff there are $k+2$ distinct elements $x_{1},\dots,
x_{k+2}$ in $\mathbb{F}_{q}$ satisfying
\begin{equation}\label{eq1}
\left\{
\begin{array}{ll}
b=x_{1}+\cdots+ x_{k+2}, \\
c=\sum_{1\leq i<j\leq k+2}x_{i}x_{j}.
\end{array}
\right.
\end{equation}
Denote $\mathbb{F}_{q}^{\ast}=\{a_{1},\ldots,a_{q-1}\}$.
\begin{itemize}
\item $b=0, c=0$.

Set $x_{i}=a_{i}, 1\leq i\leq q-1$. The conclusion holds because of the fact that all elements in $\F_q^*$ sum to 0 and Lemma \ref{lem2}.
\item $b^{2}=c,b\neq 0$.

Without loss of generality, assume that $b=-a_{1}$, and set $x_{1}=0,x_{2}=a_{2},\cdots, x_{q-1}=a_{q-1}$, then $b=-a_{1}=a_{2}+\cdots+a_{q-1}=0+x_{2}+\cdots+x_{q-1}$. From Lemma \ref{lem2},
\begin{eqnarray*}
c&=&0+b^{2}\\
&=&0+b(a_{2}+\cdots+a_{q-1})\\
&=&\big\{\sum\limits_{1\leq i< j\leq
q-1}a_{i}a_{j}\big\}-a_{1}(a_{2}+\cdots+a_{q-1}) \\
&=&\sum\limits_{2\leq
i< j\leq q-1}a_{i}a_{j}\\
&=&\sum\limits_{1\leq i< j\leq q-1}x_{i}x_{j}.
\end{eqnarray*}
Here we complete the proof that when $b^2=c$ and $k+2=q-1$, $d(u,\C)=q-k-2$.
\end{itemize}
\item $b^{2}\neq c$.

For this case, if we prove that $d(u, \mathcal{C})\neq q-k-2$ and $d(u, \mathcal{C})\leq q-k-1$, then by lemma 1 we have $d(u, \mathcal{C})= q-k-1$.
Firstly, we prove that $d(u, \mathcal{C})\neq q-k-2$.
\begin{itemize}
\item $b^{2}\neq c,b=0$.

If there are $q-1$ distinct elements $x_{1},\cdots,
x_{q-1}$ in $\F_q$ such that  $b=0=x_{1}+\cdots+x_{q-1}$. As all the nonzero elements in $\F_q$ sum to 0, then $x_{i}\neq 0,1\leq i\leq q-1$. From Lemma \ref{lem2}, $c=0$, which is contradiction with the fact that $b^{2}\neq c$.
\item $b^{2}\neq c, b\neq0$.

If there are $q-1$ distinct elements $x_{1},\cdots,
x_{q-1}$ in $\F_q$ such that  $b=x_{1}+\cdots+x_{q-1}$. As all the nonzero elements in $\F_q$ sum to 0, then there exists $x_{i}=0$. Assume that $x_{1}=0$, $x_{2}=a_{2},\cdots,x_{q-1}=a_{q-1}$, then $b=a_{2}+\cdots+a_{q-1}=-a_1$.
\begin{eqnarray*}
c&=&\sum\limits_{1\leq i< j\leq
q-1}x_{i}x_{j} \\
&=&\sum\limits_{2\leq
i< j\leq q-1}a_{i}a_{j}\\
&=&\sum\limits_{1\leq i< j\leq
q-1}a_{i}a_{j}-a_{1}(a_{2}+\cdots+a_{q-1})\\
&=&a_{1}^{2}\\
&=&b^2.
\end{eqnarray*}
Which is contradiction with $b^{2}\neq c$.
\end{itemize}
Secondly, we  prove that if $b^{2}\neq c$, then $d(u, \mathcal{C})\leq q-k-1$.
From Lemma \ref{2.3}, we know that $d(u, \mathcal{C})\leq
q-k-1$ iff there are $k+1$ distinct elements $x_{1},\dots,
x_{k+1}$ in $\mathbb{F}_{q}$ and $a\in \mathbb{F}_{q}$ satisfying
\begin{equation}
\left\{
\begin{array}{ll}
b=x_{1}+\cdots+ x_{k+1}+a, \\
c=a(x_{1}+\cdots+x_{k+1})+\sum_{1\leq i<j\leq k+1}x_{i}x_{j}.
\end{array}
\right.
\end{equation}
As $b^{2}\neq c$, then $b\neq 0$ or $c\neq 0$.  If $b\neq 0$ and $c\neq 0$, set $\eta_{1}=-b^{-1}c,\eta_{2}=0$.  If $b\neq 0$ and $c=0$, set $\eta_{1}\neq 0,-b,-2b$ and $\eta_{2}=-(b+\eta_{1})^{-1}b\eta_{1}$.  If $b=0$ and $c\neq 0$, set $\eta_{1}\neq 0, \eta_{1}^{2}\neq -c$ and $\eta_{2}=-\eta_{1}^{-1}c$. Then for each case we can check that $\eta_{1}\neq\eta_{2}$. Let $
x_{i}\in\mathbb{F}_{q}\setminus\{\eta_{1},\eta_{2}\},1\leq
i\leq q-2,a=b+\eta_{1}+\eta_{2}.$
thus, for each case, we all have that
\begin{equation*}
x_{1}+\cdots+ x_{q-2}+a=x_{1}+x_{2}+\cdots+ x_{q-2}+b+\eta_{1}+\eta_{2}
=b.
\end{equation*}
\begin{eqnarray*}
a(x_{1}+\cdots+x_{q-2})+\sum\limits_{1\leq i< j\leq
q-2}x_{i}x_{j}
&=&(b+\eta_{1}+\eta_{2})(-\eta_{1}-\eta_{2})-\eta_{1}(x_{1}+\cdots+x_{q-2})\\
&-&\eta_{2}(x_{1}+\cdots+x_{q-2}+\eta_{1})\\
&=&(-b-\eta_{1})\eta_{2}-b\eta_{1}\\
&=&c.
\end{eqnarray*}
\end{enumerate}

%\begin{thm}\label{thm2}
%Let $p=2,k+2\leq q-2$,  and $u\in \mathbb{F}_{q}^n$, $u(x)=x^{k+2}-bx^{k+1}+cx^{k}+v(x), \deg v\leq
%k-1$, then $d(u, \mathcal{C})$ satisfies the following:
%\begin{align*}
%\left\{
%\begin{array}{ll}
%2\mid k+1&\left\{\begin{array}{ll}
% 4\nmid k+1 & \left\{
%\begin{array}{ll}
%b^2\neq c & d(u, \mathcal{C})=q-k-2\\
%b^2=c,\  k+2>q/2, & d(u, \mathcal{C})=q-k-2\\
%\end{array}\right.\\
%4\mid k+1 & \left\{
%\begin{array}{ll}
% c\neq 0 & d(u, \mathcal{C})=q-k-2\\
%c=0 & \left\{
%\begin{array}{ll}
%k+2<q/2, & d(u, \mathcal{C})=q-k-2\\
%k+2\geq q/2, & d(u, \mathcal{C})\leq q-k-1
%\end{array}
%\right.
%\end{array}
%\right.
%\end{array}
%\right.\\
%2\nmid k+1 & \left\{
%\begin{array}{ll}
%4\mid k, & d(u, \mathcal{C})\leq q-k-1
%\\
%4\nmid k, &\left\{
%\begin{array}{ll}
%b\neq 0 \text{ or } c\neq 0,&  d(u, \mathcal{C})\leq q-k-1\\
%b=c=0, k+1>q/2,& d(u, \mathcal{C})\leq q-k-1
%\end{array}
%\right.
%\end{array}
%\right.
%\end{array}
%\right.
%\end{align*}
%then $d(u, \mathcal{C})\leq q-k-1$. Moreover, if $2\mid k+1$, $d(u,
%\mathcal{C})=q-k-2$.
%\end{thm}

\subsection{Proof of theorem 2 (ii).} %(details??)

\begin{enumerate}[(1)]
\item From Lemma \ref{2.3}, we know that $d(u, \mathcal{C})=q-k-2$ iff there are $k+2$ distinct elements $x_{1},\dots,
x_{k+2}$ in $\mathbb{F}_{q}$ satisfying
\begin{equation}\label{eq1}
\left\{
\begin{array}{ll}
b=x_{1}+\cdots+ x_{k+2}, \\
c=\sum_{1\leq i<j\leq k+2}x_{i}x_{j}.
\end{array}
\right.
\end{equation}
Set $x_{1}=x+y_{1},\cdots, x_{k+2}=x+y_{k+2},
m_{1}=y_{1}+\cdots+y_{k+2}, m_{2}=\sum_{1\leq i<j\leq k+2}y_{i}y_{j} $.
Then
\begin{equation}\label{eq2}
\left\{
\begin{array}{ll}
b=(k+2)x+m_{1} \\
c=\frac{(k+2)(k+1)}{2}x^{2}+(k+1)m_{1}x+m_{2}.
\end{array}
\right.
\end{equation}
Obviously, $d(u, \mathcal{C})=q-k-2$ iff the equation above has a solution.

%\item $2\mid k+1$.
\begin{itemize}
\item $2|k+1$, $4\nmid k+1$.

In this case, Equation (\ref{eq2}) can be simplified as the following equation with variable $x$.
\begin{equation}
\left\{
\begin{array}{ll}
b=x+m_{1} \\
c=x^{2}+m_{2}.
\end{array}
\right.
\end{equation}
Equation (\ref{eq2}) having a solution is equivalent to the fact that there are $k+2$ distinct elements $y_1,\cdots,
y_{k+2}$ in $\mathbb{F}_{q}$ such that $b^{2}+c=\sum_{1\leq i\leq j\leq k+2}y_{i}y_{j}.$
If $b^{2}\neq c$, as $k+2\leq q-2$ and $2|k+1$, then $k+2\leq q-3$. From Lemma \ref{lem2}
there always exist $k+2$ distinct elements $\gamma_{1},\gamma_{2},\cdots,\gamma_{k+2}$ in $\mathbb{F}_{q}^{\ast}$
such that
$$b^{2}+c=\sum_{1\leq
i\leq j\leq k+2}\gamma_{i}\gamma_{j}.$$

If $b^{2}=c$ and $k+2>q/2$, denote $t=q-1-(k+2)$. As $p=2,
2|k+1$, $4\nmid k+1$, so $4\mid
t$ and $t<q/2-1$.
From Lemma \ref{lem3}, we can get  $t$ distinct elements  $\xi_{1},\dots,
\xi_{t}$ in $\mathbb{F}_{q}^*$ such that
$$\sum_{1\leq i<j\leq t}\xi_{i}\xi_{j}=0.$$

Denote $\mathbb{F}_{q}^{\ast}=\{\xi_{1},\dots,
\xi_{t},\xi_{t+1},\dots,\xi_{q-1}\}$,
using Lemma \ref{lem2} and the fact that all the elements in $\F_q^*$ sum to 0, we can obtain
\begin{eqnarray*}
0&=&\sum\limits_{1\leq i<j\leq q-1
}\xi_{i}\xi_{j}\\
&=&\sum\limits_{1\leq i<j\leq t}\xi_{i}\xi_{j}
+(\xi_{1}+\cdots+\xi_{t})(\xi_{t+1}+\cdots+\xi_{q-1})+\sum\limits_{t+1\leq
i<j\leq q-1}\xi_{i}\xi_{j}\\
&=&0+\sum\limits_{t+1\leq i\leq j\leq q-1}\xi_{i}\xi_{j}.
\end{eqnarray*}
Therefore, in the condition that $2\mid k+1$, $4\nmid k+1$, if $b^2\neq c$ or $b^2=c$, $k+2>q/2$, then $d(u,
\mathcal{C})=q-k-2$.

%\textcolor{red}{$b^2=c$ and $k+2\leq q/2$????}

\item $4|k+1$.
We can simplify Equation (\ref{eq2}) as
\begin{equation}\label{eq3}
\left\{
\begin{array}{ll}
b=x+m_{1}; \\
c=m_{2}.
\end{array}
\right.
\end{equation}
Equation (\ref{eq3}) having a solution is equivalent to the fact that there are $k+2$ distinct elements $y_1,\cdots,
y_{k+2}$ in $\mathbb{F}_{q}$ such that $c=\sum_{1\leq i< j\leq k+2}y_{i}y_{j}.$
If $c\neq 0$, the claim holds directly from Lemma \ref{lem2}.
If $c=0$, and $k+2<q/2$, From Lemma \ref{lem3}, we can get  $k+1$ distinct elements  $\xi_{1},\dots,
\xi_{k+1}$ in $\mathbb{F}_{q}^*$ such that
$$\sum_{1\leq i<j\leq k+1}\xi_{i}\xi_{j}=0.$$
Thus, we can get $k+2$ distinct elements $0, \xi_1,\dots,\xi_{k+1}$ in $\F_q$ sum to 0.
Overall, when $4|k+1$, if $c\neq 0$ or $c=0$, $k+2<q/2$, then $d(u,
\mathcal{C})=q-k-2$.

%\textcolor{red}{$c=0, k+2\geq q/2??????$}
\end{itemize}

\item From Lemma \ref{2.3}, we know that $d(u, \mathcal{C})\leq q-k-1$ iff there are $k+1$ distinct elements $x_{1},\dots,
x_{k+1}$ in $\mathbb{F}_{q}$ and $a\in \mathbb{F}_{q}$ satisfying

\begin{equation*}
\left\{
\begin{array}{ll}
b=x_{1}+\cdots+ x_{k+1}+a \\
c=a(x_{1}+\cdots+x_{k+1})+\sum_{1\leq i<j\leq k+1}x_{i}x_{j}.
\end{array}
\right.
\end{equation*}
Denote $x_{1}=x+y_{1},\cdots, x_{k+1}=x+y_{k+1},
m_{1}=y_{1}+\cdots+y_{k+1}, m_{2}=\sum_{1\leq i<j\leq k+1}y_{i}y_{j} $.
Then,
\begin{equation*}
\left\{
\begin{array}{ll}
b=(k+1)x+a+m_{1} \\
c=\frac{(k+1)k}{2}x^{2}+km_{1}x+(k+1)ax+am_{1}+m_{2}.
\end{array}
\right.
\end{equation*}
Solving the equation system above is equivalent to solving the following equation with variable $x$.
\begin{equation}\label{eq5}
-\frac{(k+1)(k+2)}{2}x^2+x((k+1)b-(k+2)m_1)+bm_1-m_1^2+m_2-c=0.
\end{equation}
Therefore, the problem boils down to deciding if there exist $k+1$ distinct elements $y_1,\dots,
y_{k+1}$ in $\mathbb{F}_{q}$ such that the equation (\ref{eq5}) has a solution.
%\item $2\nmid k+1$.
\begin{itemize}
\item $4|k$.
In this case, Equation (\ref{eq5}) can be reduced to
\begin{equation}\label{eq6}
(x+m_1)^2+b(x+m_1)+m_2+c=0.
\end{equation}
If $ c\neq0$, we can get $k+1$ distinct elements $\{\beta_{1},\dots,
\beta_{k+1}\}\subset\F_q^*$ such that
$c=\sum_{1\leq i<j\leq k+1}\beta_{i}\beta_{j}$ according to Lemma \ref{lem2}. Then, set $y_i=\beta_i$,
$x=\beta_{1}+\cdots+\beta_{k+1}$ is  a solution of Equation (\ref{eq6}).

If $ c=0$, Denote $\alpha\in \mathbb{F}_{q}^\ast$ and $\alpha\neq b$.
According to Lemma \ref{lem2}, we can get $k+1$ distinct elements $\nu_{1},\dots,
\nu_{k+1}$ in $\F_q^*$ such that $\alpha^{2}+b\alpha=\sum_{1\leq i<j\leq
k+1}\nu_{i}\nu_{j}$. Then, set $y_{i}=\nu_{i},1\leq i\leq k+1$,  $x=\alpha+\nu_{1}+ \cdots
+\nu_{k+1}$ is  a solution of Equation (\ref{eq6}).

\item $2|k$, $4\nmid k$.
In this case, Equation (\ref{eq5}) can be reduced to
\begin{equation}\label{eq7}
b(x+m_1)+m_1^{2}+m_2+c=0.
\end{equation}
If $b\neq0$, it is easy to see that Equation (\ref{eq7})has a solution.

If $b=0, c\neq 0$, Equation (\ref{eq7})holding is equivalent to
$m_1^{2}+m_2+c=0$, which can be deduced directly for Lemma \ref{lem2}.

If $b=0, c=0$, Equation (\ref{eq7})holding is equivalent to $m_1^{2}+m_2=0$. Denote $t=q-1-(k+1)$. If  $k+1>q/2$, then $4\mid t$,and $t<q/2-1$. From Lemma \ref{lem3}, there exist
$t$ distinct $\xi_{1},\dots,
\xi_{t}$ in $\mathbb{F}_{q}^*$ satisfying
$$0=\sum_{1\leq i<j\leq t}\xi_{i}\xi_{j}.$$
Denote $\mathbb{F}_{q}^{\ast}=\{\xi_{1},\dots,
\xi_{t},\xi_{t+1},\dots,\xi_{q-1}\}$. Similarly to our proof in the first part, we can get
$\sum\limits_{t+1\leq i\leq j\leq q-1}\xi_{i}\xi_{j}=0$.

%\textcolor{red}{$b=0, c=0,k+1\leq q/2?????$}
\item $4\mid k+1$.

 In this case, equation (\ref{eq5}) can be reduced to
$$m_1x+bm_1+m_1^2+m_2+c=0.$$
Obviously, this equation has a solution.  So, if $4\mid k+1$, $d(u, \C)\leq q-k-1$.
%\item $2\mid k+1$, $4\nmid k+1$.
%
% In this case, Equation (\ref{eq5}) can be reduced to
% \begin{equation}\label{eq8}
% x^2+m_1x+bm_1+m_1^2+m_2+c=0.
% \end{equation}
%
% If $b=0, c\neq 0$, then there exist $k+1$ distinct elements $\nu_{1},\dots,
%\nu_{k+1}$ in $\F_q^*$ such that $\sum_{1\leq i\leq j\leq k+1}\nu_i=c$. Then set $y_i=\nu_i$, and $x=0$ is a solution of Equation (\ref{eq8}).
%ÏÂÃæÎÒÃÇÒªÖ¤Ã÷$\sum\limits_{t+1\leq i\leq j\leq
%q-1}\xi_{i}\xi_{j}=0$¡£ ¸ù¾ÝÒýÀí3ÖÐµÄ(iii),
%ÒýÀí4ÒÔ¼°ÓÐÏÞÓòÖÐËùÓÐÔªËØµÄºÍµÈÓÚ0Õâ¸öÊÂÊµ¿ÉµÃ,
%\begin{eqnarray*}
%0&=&\sum\limits_{1\leq i<j\leq q-1
%}\xi_{i}\xi_{j}\\
%&=&\sum\limits_{1\leq i<j\leq t}\xi_{i}\xi_{j}
%+(\xi_{1}+\cdots+\xi_{t})(\xi_{t+1}+\cdots+\xi_{q-1})+\sum\limits_{t+1\leq
%i<j\leq q-1}\xi_{i}\xi_{j}\\
%&=&0+\sum\limits_{t+1\leq i\leq j\leq q-1}\xi_{i}\xi_{j}.
%\end{eqnarray*}
\end{itemize}
\end{enumerate}
\noindent
\subsection{Proof of theorem 2 (iii) }
In order to prove the third part of Theorem 2, we have to prove the following lemmas.
\begin{lem}\label{lem4}
If $p\neq 2$,~$3<k+2< q-3$,%(??~$4\nmid q-1$ or $p\neq 2$,~$3<k+2< q-3$,~$4\mid q-1$,) ~
$k+2\neq \frac{q-1}{2}$,~$k+2\neq \frac{q+1}{2}$,
then for any $\zeta\in\mathbb{F}_{q}^{\ast}$, there exist $k+2$ distinct elements $y_{1},
\cdots, y_{k+2}$ in
$\mathbb{F}_{q}$ such that
\begin{eqnarray*}
y_{1}+\cdots+y_{k+2}&=&0\\
y_1^2+\cdots+y_{k+2}^2&=&\zeta.
\end{eqnarray*}
\end{lem}
\emph{Proof.}
\begin{enumerate}[1)]
\item $k$ is even.

If $2<k+2<\frac{q-1}{2}$, denote $y_1=\alpha, y_2=-\alpha, y_3=y, y_4=-y,
y_5=yg, y_6=-yg,\cdots, y_{k+1}=yg^{\frac{k}{2}-1},
y_{k+2}=-yg^{\frac{k}{2}-1}.$ From Lemma \ref{2.2}, the following equation with variables $\alpha, y$ has at least $q-1$ solutions in $\mathbb{F}_{q}^2$
\begin{equation}\label{eq9}
2(\alpha^2+y^2\frac{1-g^{k}}{1-g^2})-\zeta=0
\end{equation}
If $\alpha=0$ or $y=0$ or $\alpha=\pm yg^i, 0\leq i\leq
\frac{k}{2}-1$, the number of solutions in  $\mathbb{F}_{q}^2$ is at most $2k+4$. When $k+2<\frac{q-1}{2}$, $2k+4<q-1$. Then  there exists solution of equation (\ref{eq9}) satisfying  $\alpha\neq 0, y\neq 0,
\alpha\neq\pm yg^i, 0\leq i\leq \frac{k}{2}-1$. Assume that it's $(\alpha_{1},
y_{1})$. Thus, there exist $k+2$ nonzero elements $\alpha_{1},-\alpha_{1}$,~$y_{1},-y_{1}$,~$y_{1}g,
-y_{1}g \cdots, y_{1}g^{\frac{k}{2}-1}, -y_{1}g^{\frac{k}{2}-1}$ in $\mathbb{F}_{q}$ such that
\begin{equation*}
\begin{array}{ll}
\alpha_{1}+(-\alpha_{1})+y_{1}+(-y_{1})+y_{1}g+ (-y_{1}g)+
\cdots+y_{1}g^{\frac{k}{2}-1}+(-y_{1}g^{\frac{k}{2}-1})=0, \\
\alpha_{1}^{2}+(-\alpha_{1})^{2}+y_{1}^{2}+(-y_{1})^{2}+
\cdots+(y_{1}g^{\frac{k}{2}-1})^{2}+(-y_{1}g^{\frac{k}{2}-1})^{2}=
2(\alpha_{1}^2+y_{1}^2\frac{1-g^{k}}{1-g^2})=\zeta.
\end{array}
\end{equation*}
\item $k$ is odd.

According to the first part of our proof, we know that if $2<k+1<\frac{q-1}{2}$, then there exist $k+1$ distinct elements $y_{1},
\cdots, y_{k+1}$ in
$\mathbb{F}_{q}^*$ such that
\begin{eqnarray*}
y_{1}+\cdots+y_{k+1}&=&0\\
y_1^2+\cdots+y_{k+1}^2&=&\zeta.
\end{eqnarray*}
Set $y_{k+2}=0$, then the conclusion holds. So, if $k$ is odd, and $3<k+2<\frac{q+1}{2}$, the conclusion holds.
\item If $k$ is even and $\frac{q-1}{2}<k+2<q-3$, denote $t=q-k-2$, then $t$ is odd and $3<t<\frac{q+1}{2}$. From 2), we can see that if $3<t<\frac{q+1}{2}$, there exist $t$ distinct elements $z_1,\cdots, z_t$ in
$\mathbb{F}_{q}$ such that
 \begin{eqnarray*}
z_1+\cdots+z_{t}&=&0\\
z_1^2+\cdots+z_{t}^2&=&-\zeta.
\end{eqnarray*}
Note that all elements in $\F_q$ satisfy the following properties.
\begin{eqnarray*}
\sum\limits_{x\in \mathbb{F}_{q}}x&=&0 \\
\sum\limits_{x\in \mathbb{F}_{q}}x^2&=&0
\end{eqnarray*}
Then the conclusion holds.% In other words, if $k$ is even and $\frac{q-1}{2}<k+2<q-3$, the conclusion holds.

If  $k$ is odd and $\frac{q+1}{2}<k+2<q-3$, denote $t=q-k-2$, then $t$ is even and $3<t<\frac{q-1}{2}$. From 1), we can see that if $2<t<\frac{q-1}{2}$, there exist $t$ distinct elements $z_1,\cdots, z_t$ in
$\mathbb{F}_{q}$ such that
 \begin{eqnarray*}
z_1+\cdots+z_{t}&=&0\\
z_1^2+\cdots+z_{t}^2&=&-\zeta.
\end{eqnarray*}
\end{enumerate}
Similarly, we can prove the conclusion holds. %In other words, if $k$ is odd and $\frac{q+1}{2}<k+2<q-2$ the conclusion holds. (??However, when $k$ is odd and $k+2<q-2$, $k+2<q-3$.)

Hence, the proof is complete.
\begin{cor}\label{cor2}
If $p\neq 2$,~$3<k+2< q-3$, and $p\mid k+2$
then for any $\zeta\in\mathbb{F}_{q}^{\ast}$, there exist $k+2$ distinct elements $y_{1},
\cdots, y_{k+2}$ in
$\mathbb{F}_{q}$ such that
\begin{eqnarray*}
y_{1}+\cdots+y_{k+2}&=&0\\
y_1^2+\cdots+y_{k+2}^2&=&\zeta.
\end{eqnarray*}
\end{cor}
\noindent
\emph{\bf{Proof of theorem 2 (iii)}.}

%\begin{thm}(change:)
% If $p\neq 2,k+2\leq q-2$, a received word $u\in
%\mathbb{F}_{q}^n$ with corresponding polynomial
%$u(x)=x^{k+2}-bx^{k+1}+cx^{k}+v(x), \deg v\leq
%k-1$,  then
%\begin{align*}
%p\neq 2&\left\{
%\begin{array}{ll}
%p\mid k+2&\left\{
%\begin{array}{ll}
%b\neq 0 \text{ or } c\neq 0,&  d(u, \mathcal{C})= q-k-2\\
%b=c=0, k+2>q/2+1,& d(u, \mathcal{C})\leq q-k-1
%\end{array}
%\right.\\
%p\nmid k+2, &d(u, \mathcal{C})\leq q-k-1.
%\end{array}\right.
%\end{align*}
%
%
%%$d(u, \mathcal{C})\leq q-k-1$. Moreover, when $3<k+2<q-3$ and $p\mid k+2$,
%%if $b\neq 0$ or $c\neq0$, then $d(u, \mathcal{C})=q-k-2$.
%\end{thm}
%\emph{Proof.}
From Lemma \ref{2.3}, $d(u, \mathcal{C})=q-k-2$ iff there are $k+2$ distinct elements $x_{1},\cdots,
x_{k+2}$ in $\mathbb{F}_{q}$ such that
\begin{equation*}
\left\{
\begin{array}{ll}
b=x_{1}+\cdots+ x_{k+2}, \\
c=\sum_{1\leq i<j\leq k+2}x_{i}x_{j}.
\end{array}
\right.
\end{equation*}
It is equivalent to
\begin{equation}\label{eq10}
\left\{
\begin{array}{ll }
b=x_{1}+\cdots+x_{k+2},\\
b^2-2c=x_{1}^{2}+\cdots+x_{k+2}^{2}.
\end{array}
\right.
\end{equation}
Denote $x_{1}=x+y_{1},\cdots, x_{k+2}=x+y_{k+2},$ $m=y_{1}+\cdots+y_{k+2}, n=y_1^2+\cdots+y_{k+2}^2$. In order to prove  that equation (\ref{eq10}) has a solution, we just need to prove that there exist $k+2$ distinct elements $y_1,\cdots,
y_{k+2}$ in $\mathbb{F}_{q}$ such that the following equation holds.
\begin{equation}\label{eq11}
\left\{
\begin{array}{ll }
b=(k+2)x+m,\\
b^2-2c=(k+2)x^{2}+2mx+n.
\end{array}
\right.
\end{equation}
If $p\mid k+2$ and $b\neq 0$,
first according to Theorem \ref{thm1}, we can get $k+2$ distinct elements $\chi_1,\cdots,
\chi_{k+2}$ in $\mathbb{F}_{q}$ such that $b=x_1+\cdots+x_{k+2}$,
then $x=\frac{b^2-2c-n}{2b}$ is a solution of equation (\ref{eq11}).

If  $p\mid k+2$ and $b=0$, equation (\ref{eq11}) having a solution is equivalent to the fact that there are $k+2$ distinct elements $y_{1},
\cdots, y_{k+2}$ in $\mathbb{F}_{q}$ such that the following equation system holds.
\begin{eqnarray*}
y_1+\cdots+y_{k+2}&=&0\\
y_1^2+\cdots+y_{k+2}^2+2c&=&0
\end{eqnarray*}
When $c\neq0$ and $3<k+2<q-3$, it can be deduced directly from Corollary \ref{cor2}. If $k+2=3$ or $k+2=q-3$, since $p\mid k+2$, then $p=3$.
For the case $p=3$ and $k+2=3$, we know that $d(u, \mathcal{C})\leq q-k-1=q-2$ iff there are $2$ distinct elements $x_{1},
x_{2}$ in $\mathbb{F}_{q}$ and $a\in \mathbb{F}_{q}$ satisfying

\begin{equation}\label{eq3.2}
\left\{
\begin{array}{ll}
b=x_{1}+x_2+a=0 \\
c=a(x_{1}+x_2)+x_{1}x_{2}.
\end{array}
\right.
\end{equation}
As $p=3$, Equation(\ref{eq3.2}) having a solution is equivalent to $-c=(x_1-x_2)^2$ holding for distinct $x_1,x_2$. Then when $-c$ is not a nonzero square, $d(u, \mathcal{C})=q-k$.
Moreover, if $-c=(x_1-x_2)^2$ and $x_1\neq x_2$, then $a=-(x_1+x_2)\neq x_1$(otherwise $x_2=x_1$, contradiction), likewise, $a\neq x_2$. Therefore there are 3 distinct elements in $\mathbb{F}_{q}$ such that Equation (\ref{eq3.2}) holds, then $d(u,\mathcal{C})=q-k-2$. In this way, we find a new deep hole with degree $k+2$ for $p=3$ and $k=1$, which means that Cheng-Murray Conjecture doesn't hold for this special case. But $k=1$ is not usually used in designing Reed-Solomon code in practise, so this conjecture still needs further study.
For the case $p=3$ and $k+2=q-3$,  we know that $d(u, \mathcal{C})\leq q-k-1$ iff there are $k+1$ distinct elements $x_{1},\dots,
x_{k+1}$ in $\mathbb{F}_{q}$ and $a\in \mathbb{F}_{q}$ satisfying

\begin{equation*}
\left\{
\begin{array}{ll}
0=x_{1}+\cdots+ x_{k+1}+a \\
c=a(x_{1}+\cdots+x_{k+1})+\sum_{1\leq i<j\leq k+1}x_{i}x_{j}.
\end{array}
\right.
\end{equation*}
which is equivalent to $-c=\sum_{1\leq i<j\leq 4}x_ix_j$ holding for distinct $x_i\in \F_q,\ 1\leq i\leq 4$. Obviously, the equation has solutions when $q>3$.
Moreover, $d(u,\C)=q-k-2$ iff the following equation has a solution for distinct $x_1, x_2,x_3$
\begin{equation}%\label{eq3.2}
\left\{
\begin{array}{ll}
0=x_{1}+x_2+x_3 \\
c=\sum_{1\leq i\leq j\leq 3}x_ix_j
\end{array}
\right.
\end{equation}
which is equivalent to $-c=(x_1-x_2)^2$ as $p=3$. Using a similar argument, we know that $d(u,\C)=q-k-2$ iff $-c$ is a nonzero square.

From Lemma \ref{2.3}, we know that $d(u, \mathcal{C})\leq q-k-1$ iff there are $k+1$ distinct elements $x_{1},\dots,
x_{k+1}$ in $\mathbb{F}_{q}$ and $a\in \mathbb{F}_{q}$ satisfying

\begin{equation*}
\left\{
\begin{array}{ll}
b=x_{1}+\cdots+ x_{k+1}+a \\
c=a(x_{1}+\cdots+x_{k+1})+\sum_{1\leq i<j\leq k+1}x_{i}x_{j}.
\end{array}
\right.
\end{equation*}
Denote $k+1 \equiv r~\mod~p, x_{1}=x+y_{1},\cdots, x_{k+1}=x+y_{k+1}, m=y_{1}+\cdots+y_{k+1},
n=y_1^2+\cdots+y_{k+1}^2$.
Then,
\begin{equation*}
\left\{
\begin{array}{ll}
b=(k+1)x+a+m\\
c=\frac{(k+1)k}{2}x^{2}+kmx+(k+1)ax+am+\frac{m^2-n}{2}.
\end{array}
\right.
\end{equation*}
Solving the equation system above is equivalent to solving the following equation with variable $x$.
\begin{equation}\label{eq12}
(r+r^{2})x^{2}+x(2m+2mr-2br)+m^{2}+n-2bm+2c=0
\end{equation}
We discuss the solution of Equation (\ref{eq12}) according to the following three cases.
\begin{enumerate}[(i)]
\item $r=0$.

In this case, Equation (\ref{eq12}) can be reduced to the form $2mx+m^2+n-2bm+2c=0$. Obviously, it has a solution.
\item $r=-1$.

In this case, $p\mid k+2$, thus, for  $b\neq 0$ or $b=0,\ c\neq 0, \ p\neq 3$ or $p=3$, we have discussed due to the first part of our proof. For $b=0$ and $c=0$, denote $t=q-1-(k+1)$, then $p|t$. From Corollary \ref{cor1}, if $k+1>q/2$, there are $k+1$ distinct elements $\xi_{1},\dots,
\xi_{k+1}$ in $\mathbb{F}_{q}^*$ such that $$\sum_{1\leq i\leq j\leq k+1}\xi_{i}\xi_{j}=0.$$
Then Equation (\ref{eq12}) has a solution.

%\textcolor{red}{$b=0, c=0, k+2<q/2+1??????$}
\item $r\neq -1, 0.$

Denote $S=\{x^2~|~x\in\mathbb{F}_{q}\}$. As we know, equation $ax^2+bx+c=0 (a\neq
0)$ over finite field $\mathbb{F}_{q}$ with characteristic $p\neq
2$ having a solution is equivalent to discriminant $D=b^2-4ac\in S.$
For Equation (\ref{eq12}), The determinant
\begin{eqnarray*}
D_{1}&=&(2m+2mr-2br)^{2}-4(r+r^{2})(m^{2}+n-2bm+2c)\\
&=&4[(1+r)(m^{2}-rn-2cr)+b^{2}r^{2}].
\end{eqnarray*}

If $2<k+1<\frac{q+1}{2}$ and $k+1$ is even, denote
$$A_{1}=\big\{\frac{\alpha^{2}}{1+r}-\frac{b^{2}r^{2}}{1+r}+2cr~|~ \alpha \in
\mathbb{F}_{q}\big\}.$$
According to Lemma \ref{4.3}, we can obtain that there are $k+1$ distinct nonzero elements $y_{1},
\cdots, y_{k+1}$ in $\mathbb{F}_{q}$ satisfying
$$(m^{2}-rn)\in A_{1},$$
which is equivalent to $D_{1}=4[(1+r)(m^{2}-rn-2cr)+b^{2}r^{2}]\in S$.

Similarly, if $k+1$ is odd, we can get $k$ distinct nonzero elements $y_{1},
\cdots, y_{k}$ in $\mathbb{F}_{q}$ satisfying
$(m^{2}-rn)\in A_{1},$ then set
$y_{k+1}=0$.

If $q-3\geq k+1\geq\frac{q+1}{2}$, denote $t=q-k-1, t \equiv s\mod p$,
then $s\neq 0, 1$ and $1<t\leq \frac{q-1}{2}<\frac{q+1}{2}$.
Set $m_{1}=y_{1}+\cdots+y_{t}, n_{1}=y_1^2+\cdots+y_{t}^2$,
when $s\neq 0, 1, 1<t<\frac{q+1}{2}$, the determinant of the following equation with variable $x$
\begin{equation}
(s-s^{2})x^{2}+x(2m_1-2m_1s-2bs)-m_1^{2}+n_1-2bm_1-2c=0
\end{equation}
is
\begin{eqnarray*}
D_{2}&=&(2m_1-2m_1s-2bs)^{2}-4(s-s^{2})(-m_1^{2}+n_1-2bm_1-2c)\\
&=&4[(1-s)(m_1^{2}-sn_1+2cs)+b^{2}s^{2}].
\end{eqnarray*}
From Lemma \ref{4.3}
there exist $t$ distinct elements in  $\mathbb{F}_{q}$ such that
$$D_{2}\in S.$$
Then Equation (16) has a solution in $\mathbb{F}_{q}$, which is equivalent to the fact than the following equation has a solution in
$\mathbb{F}_{q}^{t+1}$.
\begin{equation}
\left\{
\begin{array}{ll }
-b=x_{1}+\cdots+x_{t}-a,\\
-b^{2}+2c=x_{1}^{2}+\cdots+x_{t}^{2}-a^{2}.
\end{array}
\right.
\end{equation}
Our conclusion holds directly by the following property of finite fields. When $q>3$,
\begin{equation*}
\left\{
\begin{array}{ll }
\sum\limits_{x\in \mathbb{F}_{q}}x=0 \\
\sum\limits_{x\in \mathbb{F}_{q}}x^2=0
\end{array}
\right.
\end{equation*}
%¼´¿ÉÖ¤µÃ¡£ \qed
%\section{ÖÂÐ»}
\end{enumerate}
%\section{Conclusion}
\emph{Acknowledgement.} We express our deep gratitude to professor Daqing Wan for his useful suggestions. And we also acknowledge with gratitude the support by National Development Foundation for Cryptological Research(No. MMJJ201401003).

\end{document}